\documentclass{amsart}
\usepackage[utf8]{inputenc}

\usepackage[english]{babel}

\usepackage{amssymb}
\usepackage{amsmath}
\usepackage{amscd}
\usepackage{amsthm}
\usepackage{enumerate}
\usepackage{cite}
\usepackage{xcolor}
\usepackage{graphicx}

\newcommand{\Pj}{\mathbb{P}}

\newcommand{\T}{\mathbb{T}}

\newcommand{\Oc}{\mathcal{O}}

\newcommand{\Ic}{\mathcal{I}}

\newcommand{\Sc}{\mathcal{S}}

\DeclareMathOperator{\red}{red}

\DeclareMathOperator{\reg}{reg}

\newtheorem{theorem}{Theorem}[section]
\newtheorem{proposition}[theorem]{Proposition}
\newtheorem{lemma}[theorem]{Lemma}
\newtheorem{corollary}[theorem]{Corollary}

\theoremstyle{definition}
\newtheorem{definition}[theorem]{Definition}
\newtheorem{notation}[theorem]{Notation}
\newtheorem{example}[theorem]{Example}
\newtheorem{remark}[theorem]{Remark}
\begin{document}

\title{Terracini loci of curves}

\author[E.~Ballico]{Edoardo Ballico}
\address{E. Ballico, Dipartimento di Matematica, Universit\`a di Trento, Italy}
\email{edoardo.ballico@unitn.it}
\author[L.~Chiantini]{Luca Chiantini}
\address{L. Chiantini, Dipartimento di Ingegneria dell'Informazione e Scienze Matematiche, Universit\`a di Siena, Italy}
\email{luca.chiantini@unisi.it}

%\author{Edoardo Ballico         \and Luca Chiantini
   %     }


%\institute{E. Ballico \at
%              Dipartimento di Matematica, Universit\`a di Trento, Italy
  %            \email{edoardo.ballico@unitn.it}           %  \\
 %          L. Chiantini \at Dipartimento di Ingegneria dell'Informazione e Scienze Matematiche, Universit\`a di Siena, Italy 
  %           \email{luca.chiantini@unisi.it}
%}


\maketitle

\begin{abstract}
We study subsets $S$ of curves $X$ whose double structure does not impose independent conditions to a linear series $L$,
but there are divisors $D\in |L|$ singular at all points of  S.
These subsets form the Terracini loci of $X$. We investigate Terracini loci, with a special look  towards their non-emptiness, mainly
in the case of canonical curves, and in the case of space curves.
\end{abstract}

\section{Introduction}
\label{intro}
Terracini loci $\T(X,L,x)$ of a projective variety $X$ (over an algebraically closed field of characteristic $0$)
are subsets of  the set $\Sc(X,x)$ of all reduced finite subsets 
$S\subset X_{\reg}$ of cardinality  $x$, with the property that the double scheme $2S$ on $S$ does not impose independent conditions to the linear series $L$,
and there are divisors in $|L|$ passing through $2S$.
Terracini loci are certainly involved in the study of interpolation properties of the image of $X$ in the map induced by $L$, but they 
also assume great importance in the theory of secant varieties to embedded varieties $X$, for they  are connected with points of
the abstract secant variety of $X$ in which the differential of the map to the embedded secant variety drops rank.
More details on the initial properties of Terracini loci can be found in \cite{BallC}.

In this paper, we focus the attention to Terracini loci of curves. Even if curves are never defective, so that their secant varieties always have the expected dimension,
yet there are special points $p$ in which the  Terracini Lemma fails to provide the dimension of the Zariski tangent space to the secant
variety. This happens typically when the set of points in $X$ which generates $p$ belongs to some Terracini locus. 

We are mainly concerned with the problem of the non-emptiness of a Terracini locus $\T(X,L,x)$. Certainly $\T(X,L,x)$ is empty when
$X$ is a rational normal curve in $\Pj^r$, $L$ is the complete  hyperplane linear series and $x\leq (r+1)/2$, because every subscheme
of length $\leq r+1$ in $X$ is linearly independent. We will see that, for odd $r$, the emptiness of $\T(X,L,(r+1)/2)$ characterizes rational
normal curves (Proposition \ref{ratiomax}). This means that there is a large mass of examples of non-empty Terracini loci of curves.

Our analysis considers mainly two cases: canonically embedded curves (Section 4) and curves in $\Pj^3$ (Section 5). 

In the former case, the 
series $L$ is the complete canonical linear series $K_X$, and one can easily expect that the  existence of sets in   $\T(X,K_X,x)$ strongly depends
on the geometry of linear series on the curve. We show indeed that subsets of theta-characteristics on $X$ prove that $\T(X,K_X,x)$
is non-empty when $g/2\leq x\leq g-1$. For $x<g/2$ we see (Proposition \ref{gon}) that the non-emptiness of $\T(X,K_X,x)$ is linked to the
gonality of $X$.

For space curves, we use induction on the degree and genus, with the technique of smoothing reducible curves, to prove that  $\T(X,L,x)$
is non-empty for some smooth curves of degree $d$ and genus $g\leq d-3$, when $6\leq 2x\leq d$  (Theorem \ref{i1}). We also prove the
existence of smooth special curves with non-empty Terracini loci in some components of the Hilbert scheme of  space curves  even outside the Brill-Noether range
(Theorem \ref{i1.1} and Remark \ref{ll0}). 

We would like to thank the anonymous referees for precious observations on the preliminary versions of the paper.

\section{Preliminaries}\label{sec:notation}

We work over an algebraically closed field of characteristic $0$.

\subsection{Notation}

For any $0$-dimensional subscheme $S$ we denote with $\ell(S)$ the length of $S$. When $S$ is reduced,
then $\ell(S)$ is the cardinality of $S$.
\smallskip

Let $X$ be an integral projective variety. For any point $P\in X_{\reg}$, with homogeneous maximal ideal $m_{P,X}$, 
we denote with $(2P,X)$ the subscheme of $X$ defined by $m_{P,X}^2$. 
If $S=\{P_1,\dots,P_x\}\subset X$ is a finite set of points, then we denote with $(2S,X)$ the non-reduced scheme 
$2S=\bigcup_{i=1}^x (2P_i,X)$.\\
When it is clear which $X$ we refer to, then we will write simply $2S$ instead of $(2S,X)$.
\smallskip

For any scheme $X$ and for any integer $x$ we denote with $\Sc(X,x)$ the set of all reduced finite subsets 
$S\subset X_{\reg}$ of length $x$. $\Sc(X,x)$ is an open subset of the symmetric product $X^{(x)}$.
\smallskip

Let $X$ be a projective curve, and let $S$ be a finite subset of $X_{\reg}$. When $S$ identifies a Cartier divisor on $X$, by abuse,
we will continue to denote it with $S$. We will also denote with $|S|$ the complete linear series associated to $S$, and we will denote
with $h^0(S)$ the dimension of the space of sections of the associated line bundle. We will use the same convention for $2S$.

For any scheme $Z\subset \Pj^r$ let $\langle Z\rangle \subseteq \Pj^r$ denote the minimal linear subspace of $\Pj^r$ containing $Z$.

\subsection{Terracini loci}

We recall from \cite{BallC} the definition of Terracini locus of a projective variety. 

\begin{definition}\label{defT}
Let $X$ be an integral projective variety, $L$ a line bundle on $X$ and $V\subseteq H^0(L)$ a linear subspace. 
Set $m:= \dim X$. Fix $S\in \Sc(X,x)$. 
We say that $S$ is in the Terracini locus $\T(X,L,V,x)$ if $V(-(2S,X))\ne (0)$ and 
$\dim V(-(2S,X)) > \dim V -x(m+1)$.\\
We say that the integer $\delta (S,X,L,V):= \dim V(-(2S,X)) -\dim V +x(m+1)$ is the (Terracini) defect of $S$ 
with respect to $(L,V)$. \\
When $V=H^0(L)$, we will drop $V$ in the notation.
\end{definition}

We will consider, throughout the paper, mainly the case where $X$ is a curve, i.e. $m=1$. In this situation a finite set $S$ of length 
$x$ lies in the Terracini locus  when $\dim V(-(2S,X)) > \max\{\dim V -2x,0\}$.

We can extend the definition of Terracini loci to include some non-reduced 0-dimensional subschemes.

\begin{definition}\label{defTT}
Let $C\subset \Pj^r$ be a smooth and connected non-degenerate curve. For each positive integer $x$ let $C^{(x)}$ denote the symmetric product of $x$ copies of $C$. The variety $C^{(x)}$ is a connected projective variety of dimension $x$ parametrizing the degree $x$ zero-dimensional schemes of $X$, and  there is a non-empty Zariski open subset $S(C,x) \subset C^{(x)}$ which parametrizes subsets of cardinality $x$. For each positive integer $x$ let $\tilde{\T}(C,x)$ denote the
set of all $Z\in C^{(x)}$ such that $\dim \langle 2Z\rangle \le 2x-2$ and $\langle 2Z\rangle \ne \emptyset$. Obviously, $\T(C,1) =\tilde{\T}(C,1)=\emptyset$. The inclusion $S(C,x)\subset C^{(x)}$ induces an inclusion $\T(C,x)\subseteq \tilde{\T}(C,x)$.
The semicontinuity theorem for cohomology gives that $\tilde{\T}(C,x)$ contains the closure of $\T(C,x)$ in $C^{(x)}$. Sometimes, $\tilde{\T}(C,x)$ is not equal to the closure of $\T(C,x)$ in $C^{(x)}$ (Example \ref{exe1}).
\end{definition}

\begin{example}\label{exe1}
Fix an even integer $d\ge 4$, a line $D\subset \Pj^2$ and $p\in D$. There is a smooth degree $d$ curve $C\subset \Pj^2$ such that $p\in C$ and $D\cap C = dp$, i.e. $p$ is a total ramification point of $C$ with $D$ as its tangent line.
Since $C$ has only finitely many multitangent lines, $\T(C,x)$ is finite and hence it is closed in $\tilde{\T}(C,x)$. Since $d$ is even and $d\ge 4$, $\frac d2 p$ lies in $ \tilde{\T}(C,d/2)\setminus \T(C,d/2)$.
\end{example}

\section{Generalities}\label{sec:gen}
\begin{example} Let $X\subset \Pj^N$ be an integral curve. Here we take $L=\Oc_X(1)$ and $V=$ the linear series 
of  hyperplanes. For a general $S=\{P_1,\dots,P_x\}\subset X_{\reg}$ write $2S$ for $(2S,X)$.
The elements of $V(-2S)$ correspond to hyperplanes containing the tangent lines to $X$ at the points $P_i$'s. 
It follows that $S$ lies in the Terracini locus
when the span of the $x$ tangent lines at the points of $X$  has dimension smaller than the expected one.
Then by definition $\T(X,L,V,1)$ is always empty.
\end{example}

On the other hand, when the map induced by the linear series $V$ is not birational onto the image  the Terracini locus
$\T(X,L,V,1)$ may contain some points.

\begin{example}\label{hypell} Let $X$ be a smooth hyperelliptic curve of genus $g\ge 2$. We can describe all the Terracini loci with respect 
to $L=K_X$ and  $V=H^0(K_X)$. \\
Let $B\subset X$ be the set of the Weierstrass points of $X$. For any integer $x$ let $\Sc(B,x)$ denote the set of all subsets of $B$ with cardinality $x$. The set $B$ is the ramification locus of the $g^1_2$ of $X$. Since we work in characteristic $0$,  $\ell(B)=2g+2$.  
Fix a positive integer $x$ and take $S\in S(X,x)$. If $x\ge g$, then $\deg (K_X(-(2S,X))) <0$ and hence $h^0(K_X(-(2S,X))) =0$. Thus 
$\T(X,K_X,x)=\emptyset$ for all $x\ge g$. \\
Assume that $1\le x\le g-1$. Let $h: X\to \mathbb{P}^1$ be the morphism associated to the $g^1_2$ of $X$.
The linear system $|K_X|$ is the minimal sum of $g-1$ copies of the $g^1_2$ of $X$. 
Hence every base-point free special line bundle on $X$ is the sum of at most $ g-1$ copies of the $g^1_2$. 
Thus $\T(X,K_X,1)=B$ and $\T(X,K_X,g-1)) = \Sc(B,g-1)$. More generally, $S$ belongs to $\T(X,K_X,x)$ if
 and only if we have that either $S\cap B\ne \emptyset$ or there are $p, q\in S$ such that $p\ne q$ and $h(p) =h(q)$ for
these conditions are equivalent to $h^0(K_X(-(2S,X))) >g-2x$. \\
It is easy to realize that, for each $x\in \{2,\dots ,g-2\}$, the elements of $\T(X,K_X,x)$ with (maximal) defect 
$x$ are the elements of $\Sc(B,x)$. If $e:= \#(S\cap B)$ and $S$ has $f$ distinct sets $\{p_i,q_i\}$, $1\le i\le f$, with $p_i\ne q_i$ and $h(p_i)=h(q_i)$, then $h^0(K_X(-2S)) =g-2x+e+f$.
\end{example}

For the rest of the section let us go back to the case in which  $X\subset \Pj^N$ is an integral curve,
and  we take $L=\Oc_X(1)$ and $V=$ the linear series of  hyperplanes.

\begin{example}
Let us consider what happens when $x=2$
and  $X$ is a plane curve. Thus $\max\{\dim V -2x,0\}=0$. Then $S=\{P,Q\}\in\Sc(2)$ belongs to the Terracini locus if and
only if the tangent lines to $X$ in $P$ and $Q$ coincide. Since in characteristic $0$ not every tangent line is bitangent,
then $\T(X,L,V,2)$ is either empty or finite.\\
If $N>2$, the tangent lines to two general points of $X$ span a $3$-dimensional linear subspace 
(recall that we work in characteristic $0$). The set $S=\{P,Q\}$ lies in the Terracini locus $\T(X,L,V,2)$ when the 
tangent lines in $P,Q$ meet at some point $P_0$, i.e. there exists a plane containing the two tangent lines. 
In this case, the projection of $X$ from $P_0$ is a curve with (at least) two cusps. 
\end{example}

\section{Canonically embedded curves}\label{sec:can}
Let us turn now to the case where $X$ is a smooth curve of genus $g$ and we consider the complete canonical linear series $L=K_X$. 
We will describe in several cases the locus $\T(X,K_X,x)$.\\
Since we already treated the case of hyperelliptic curves in Example \ref{hypell}, we assume that $g\geq 3$ and $K_X$ embeds $X$ in $\Pj^{g-1}$.
\smallskip

We start with a very easy observation, which shows that we need to distinguish two cases, depending if $x$ is smaller than $g/2$ or not.

\begin{proposition}\label{RR} A reduced set $S$ of length $x$ lies in $\T(X,K_X,x)$ if and only if either $2x< g$ and $h^0(2S)>1$, i.e. the linear series $|2S|$
is not a singleton, or $2x\geq g$ and $h^0(K_X-2S)>0$, i.e. $2S$ is special. 
\end{proposition} 
\begin{proof} By definition the set $S$  lies in $\T(X,K_X,x)$ if and only if $h^0(K_X-2S)>\min\{0, h^0(K_X)-2x\}$.  Since $h^0(K_X)=g$, we distinguish between 
 $2x<g$ and  $2x\geq g$. In the latter case  $S$  lies in $\T(X,K_X,x)$ if and only if $h^0(K_X-2S)>0$. In the former case, since by Riemann-Roch
 $h^0(2S)=2x-g+1+h^0(K_X-2S)$, $S$  lies in $\T(X,K_X,x)$ if and only if $h^0(2S)>1$. 
\end{proof}

It follows from the previous proposition that the Terracini locus $\T(X,K_X,x)$ is empty if $x\geq g$, because in this case the degree of $K_X-2S$ is
negative.\\
Let us consider the extremal case $x=g-1$.

\begin{example}\label{g-1}
A set $S$ (resp. scheme)  of length $g-1$ belongs to  $\T(X,K_X,g-1)$ (resp. $\tilde{\T}(X,K_X,g-1)$) if and only if $h^0(K_X-2S)>0$ which, for degree reasons, implies that $2S$ is a canonical divisor. Thus subsets $S\in \T(X,K_X,g-1)$ correspond to divisors in some non empty linear series $G$ such that $2G=K_X$, i.e. a \emph{theta-characteristic} of $X$. \\
It is well known (\cite{Harr82}) that $X$ has a finite number, exactly $2^{2g}$, theta-characteristics. A theta-characteristic is odd or even, depending on the parity of $h^0(G)$.
The number of odd theta-characteristics is $2^{g-1}(2^g-1)$ while there are $2^{g-1}(2^g+1)$ even theta-characteristics.\\
Now assume that $X$ is general in the moduli space $\Mc_g$. In this case,  by \cite{Harr82} Corollary 1.11,
$h^0(G)\le 1$ for every theta-characteristic $G$ on $X$, and for each odd
theta-characteristic $G$ on $X$ the divisor $D$ with $\{D\} =|G|$ is reduced. Thus for  $X\in \Mc_g$ general the Terracini locus $\T(X,K_X,g-1)$
is finite, of cardinality $2^{g-1}(2^g-1)$. \\
There are $X\in \Mc_g$ with theta-characteristics $G$ such that $h^0(G)\geq 2$. For such curves $\T(X,K_X,g-1)$ is infinite.\\
On the contrary, a natural question is to ask if there are $X\in \Mc_g$ such that $\T(X,K_X,g-1)=\emptyset$, i.e. no reduced divisor is the zero-locus of an
effective theta-characteristic.\\ 
In the case $g=3$ this is equivalent to ask if there is a smooth degree $4$ plane curve
$X$ with $28$ flexes of higher order, i.e. $28$ lines $L\subset \Pj^2$ meeting $X$ at a unique point. The total
weight of all flexes of a smooth plane quartic is $24$, because its Hessian determinant has degree $6$. Thus there is no such
$X$ for $g=3$. For any $X$ we have $\#\tilde{\T}(X,K_X,g-1)\ge 2^{g-1}(2^g-1)$ and either $\#\tilde{\T}(X,K_X,g-1)= 2^{g-1}(2^g-1)$ (case $h^0(G)\le 1$ for all theta-characteristic $G$ of $X$) or $\dim \tilde{\T}(X,K_X,g-1)>0$.
\end{example}

In any case, taking subsets of a theta-characteristic, we obtain immediately the following

\begin{proposition}\label{lessg-1}
Fix an integer $x$ such that $g/2 \leq x\leq g-1$. Then the Terracini locus $\tilde\T(X,K_X,x)$ is non empty.
\end{proposition}

Of course when $X$ has a positive dimensional theta-characteristic, then also $\T(X,K_X,x)$ is infinite for all $x$ between $g/2$ and $g-1$.

\begin{remark}
Fix an integer $x$ such that $g/2\le x\le g-2$. One can try to extend elements of $\T(X,K_X,x)$ to elements  $\T(X,K_X,x+1)$, with the addition
of suitable points. Respect to this, we can observe:\\ 
\quad (a)  Assume that there exists $S\in \T(X,K_X,x)$ with $\dim \langle 2S\rangle = g-2$, i.e. $H:= \langle 2S\rangle$ is a hyperplane. The point $p\in X\setminus S$
has the property that $S\cup \{p\}$ lies in $\T(X,K_X,x+1)$ if and only if $p\in H$ and $H$ is tangent to $X$ at $p$. 
There is $Z\in \tilde{\T}(K_X,x+1)$ containing $S$
if and only if either $H$ is tangent to $X$ at a point of $X\setminus S$ or there is $o\in S$ such that the scheme $H\cap X$ contains $o$ with multiplicity at least $4$.\\
\quad (b) Assume that there exists $S\in \T(X,K_X,x)$ with $\dim \langle 2S\rangle = g-3$. The differential of the rational map $\phi: X\setminus X\cap \langle
2S\rangle \to \Pj^1$ induced by the linear projection from $\langle 2S\rangle$ shows that there are only finitely many
$p\in X\setminus S$ such that $S\cup \{p\}\in \T(X,K_X,x+1)$. There are at least $2$ such points $p$, because $\phi$ extends
to a morphism $\psi : X\to \Pj^1$ by the smoothness of $X$ and $X$ has at least $2$ ramification points, because $g>0$.\\
\quad (c) Assume that there exists $S\in \T(X,K_X,x)$ with $\dim \langle 2S\rangle =g-4-2k-\epsilon$, with $k$ non-negative integer
 and  $\epsilon\in \{0,1\}$. Then $S\cup \{p_1,\dots ,p_{k+1}\}$ lies in $\T(X,K_X,x+k+1)$ for a general
$(p_1,\dots ,p_{k+1})\in X^{k+1}$.
\end{remark}

Let us now consider the case $x<g/2$.

\begin{proposition} Assume $x<(g-1)/2$, and assume that $\T(X,K_X,x)$ is non empty. Then
$$  \dim (\T(X,K_X,x+1))\geq 1+\dim( \T(X,K_X,x)).$$
\end{proposition}
\begin{proof} Pick $S\in \T(X,K_X,x)$ and $p$ general in $X$. Since $2x< g$ then by Proposition \ref{RR} we
have $h^0(2S)\geq 2$. It follows $h^0(2(S\cup\{p\}))\geq h^0(2S)\geq 2$. \\
When $x+1<g/2$ , the last inequality is sufficient to conclude that $S\cup\{p\}$ sits in $T(X,K_X,x+1)$, by Proposition \ref{RR} again,
thus the inequality on the dimensions holds.\\
Assume $x=(g-2)/2$. Since $h^0(2(S\cup\{p\}))\geq 2$, then by Riemann-Roch
$$h^1(2(S\cup\{p\})) \geq g-1-2\bigg(\frac{g-2}2 +1\bigg) +2 =1.$$
Thus $2(S\cup\{p\})$ is special, hence it belongs to $T(X,K_X,x+1)$, by Proposition \ref{RR}, and we conclude as before.
\end{proof}

\begin{corollary} \label{nobuchi}
Fix the minimal integer $x$ such that  the Terracini locus $\T(X,K_X,x)$ is non empty. Then for all $y$ with $x\leq y\leq g-1$ the
Terracini locus $\T(X,K_X,y)$ is non-empty.
\end{corollary}

We saw in the previous section that $\T(X,K_X,1)$ is non-empty if and only if $X$ is hyperelliptic. 
The last $2$ propositions of this section link the gonality $k$ of $X$ with the minimal $x$ such that $\T(X,K_X,x)\neq \emptyset.$

\begin{proposition}\label{gon} If $\T(X,K_X,x)\neq \emptyset$, then $X$ has a linear series $g^1_{2x}$. 
For the converse,  if $X$ has a linear series $g^1_k$ with $k< g/2$, then  $\dim(\T(X,K_X,k))\geq 1$.
\end{proposition}
\begin{proof} Recall that we are assuming $g\geq 3$ and $X$ not hyperelliptic.\\
The first assertion is trivial if $x\geq g/2$, for every $X$ satisfies  $\T(X,K_X,x)\neq \emptyset$ and every $X$ has a linear series $g^1_{2x}$.
When $x<g/2$, the first assertion follows immediately by Proposition \ref{RR}.\\
For the second assertion, consider the degree $k$ map $f: X\to \Pj^1$ associated to the linear series $g^1_k$. 
For any $p\in \Pj^1$ let $D_p:= f^{-1}(p)$ denote the associated degree $k$ divisor. Note that $h^0(2D_p)\ge
h^0(\Oc_{\Pj^1}(2p)) =3$. Then, as in Proposition \ref{RR}, a general $D_p$ belongs to $\tilde\T(X,K_X,k)$.
Since we work in characteristic $0$, then a general $D_p$ is formed by $k$ distinct points, hence it belongs to $\T(X,K_X,k)$.
\end{proof}

Before we can refine the previous proposition, let us see what happens for trigonal curves.

\begin{example}\label{trigonal}
Let $X$ be a smooth trigonal curve of genus $g\ge 4$, canonically embedded in $\Pj^{g-1}$,
 and let $f: X\to \Pj^1$ be the degree $3$ morphism associated to the $g^1_3$. 
By the Castelnuovo's inequality \cite{Kani84}, $f$ is unique if $g\ge 5$. Let $\Sigma\subset X$ denote the set of all ramification points of the map
and let $\Sigma'\subseteq \Sigma$ denote the set of all $p\in \Sigma$ which belong to fibers of cardinality $2$. 
The points $p\in\Sigma \setminus \Sigma'$ are called the \emph{total ramification points} of $f$, because the fiber containing
$p$ is supported at $p$.\\ 
$\Sigma'=\Sigma$ if $X$ is a general trigonal curve of genus $g$, but there are trigonal curves in which the equality
fails, and also trigonal curves with $\Sigma' =\emptyset$, e.g. the degree $3$ Galois coverings of $\Pj^1$. \\
Take $p\in \Sigma'$ and consider the point $q_p\neq p$ in the fiber of $f$ through $p$. We claim that
$S=\{p,q_p\}$ belongs to $\T(X,K_X,2)$, so that $\T(X,K_X,2)$ is non-empty.\\
To prove the claim, consider that $h^0(2p+q_p)=2$. Then  $h^0(2p+2q_p)=h^0(2S)\geq 2$, which proves the claim when
$g>4$, by Proposition \ref{RR}. For $g=4$, since $h^0(2S)\geq 2$, the divisor $2S$ is special by Riemann-Roch, and the 
claim follows again by Proposition \ref{RR}.
\end{example}

\begin{definition} Given a linear series $G$ on $X$ which is a $g^1_d$, $G$ is \emph{tamely ramified}
if there is a non-reduced divisor $D$ in $G$ in which all the points appear with coefficient $\leq 2$.
\end{definition}

\begin{proposition} Let $X$ be a smooth  curve of genus $g\ge 4$, canonically embedded in $\Pj^{g-1}$.
Assume that $X$ has a base point free pencil $R\in \mathrm{Pic}^k(X)$ such that $h^1(R^{\otimes 2}) >0$,  and let $f: X\to \Pj^1$ be the degree $k$ morphism induced by $|R|$. 
Then $\tilde{\T}(X,K_X,k-1)\ne \emptyset$.
Assume that the $f$ is tamely ramified. Then $\T(X,K_X,k-1)\ne \emptyset$.
\end{proposition}
\begin{proof} If $k-1\geq g/2$ the claim follows from  Proposition \ref{lessg-1}. 

Assume $k-1<g/2$. Since $g>0$ and $\Pj^1$ is algebraically simply connected,
there is a divisor $A = m_1q_1+\cdots +m_sq_s\in |R|$ with $m_1\ge 2$. Since $2A$ is a special divisor, $A\in \tilde{\T}(X,K_X,k-1)$.
Now assume that $f$ is tamely ramified. Thus all $A\in |R|$ have profiles with multiplicities $1$ or $2$.
Since $g>0$, $A=2p_1+\dots+2p_i+p_{i+1}+\dots + p_j\in |R|$ with $p_1,\dots ,p_j$  distinct points 
and $i\geq 1$. Note that $i+j=k$, so that $j<k$. Then $S=\{p_1,\dots,p_j\}$ sits in $\T(X,K_X,j)$.
Namely $j\leq k-1<g/2$ and $h^0(R)\geq 2$ implies $h^0(2S)\geq 2$. Then $\T(X,K_X,j)\ne \emptyset$, and the claim follows from Corollary \ref{nobuchi}.
\end{proof}

\section{Embedded curves}\label{sec:emb}

In this section we consider reduced and locally complete intersection curves $X$ embedded in a projective space $X\subset \Pj^r$.

\begin{notation}
We denote with $N_X$ the normal bundle of $X$.\\
For brevity, we will denote with $\T(X,x)$ the Terracini locus of $X$ with respect to the (non-necessarily complete) linear series 
of hyperplane divisors.\\
For all integers $g\ge 0$, $r\ge 3$ and $d\ge r+g$, we will denote with $H(d,g,r)$  the set of all smooth and non-degenerate curves
$X\subset \Pj^r$ of degree $d$ and genus $g$ \emph{such that $h^1(\Oc_X(1))=0$}. \\
Fix $X\in H(d,g,r)$. Since $X$ is smooth, $N_X$ is a quotient of $O_X(1)^{\oplus r}$ by the Euler's sequence. 
Thus $h^1(N_X) =0$. Therefore $H(d,g,r)$ is smooth, and
$$\dim H(d,g,r) = (r+1)d+(r-3)(1-g).$$
Since all  non-special line bundles of degree $d$ on a curve of genus g have the same number of sections,
and all smooth genus g curves are parametrized by an irreducible variety, $H(d,g,r)$ is irreducible.
See \cite{Ein87}.\\
An example of J. Harris, generalized by L. Ein in   \cite{Ein87}, shows that both claims on the  irreducibility and the 
dimension fail, for large $r$, if we drop the assumption on the vanishing of $h^1(\Oc_X(1))$. On the other hand, 
we will work mainly with curves in $\Pj^3$, and in this case we can drop the assumption by \cite{Ein86}. 
\end{notation}

It is immediate, by Bezout formula, that $\T(X,x)$ is empty if $x>d/2$. So we analyze the case $2x\le d$.\\

\begin{remark}\label{plus1} 
Fix $x$ with $2x<r$ and let $X\subset\Pj^r$ be a non-degenerate irreducible curve such that $\T(X,x)\neq \emptyset$. Then also
$\T(X,x+1)\neq \emptyset$.\\
Indeed if $S\in\T(X,x)$, then the  scheme $(2S,X)$ lies in a linear system of hyperplanes of dimension at least $r-2x+1\geq 2$. Thus, for $p\in X $ general,
the  scheme $(2(S\cup\{p\}),X)$ lies in a linear system of hyperplanes of dimension at least $r-2x-1\geq 0$.
\end{remark}

\begin{notation} Take $X\subset \Pj^r$. An \emph{arrow} in $\Pj^r$ is a non-reduced scheme of length $2$. The set of all arrows in $\Pj^r$
supported at $p$ is closed in the Hilbert scheme, and it has dimension $r-1$. Thus the set of all arrows in $\Pj^r$ has dimension $2r-1$.\\
Note that an arrow $w$ with support  $p$ determines a line $r_w$ through $p$. If $p\in X$ and $X$ is  smooth at $p$, then $r_w$ is tangent
to $X$ exactly when $X$ contains the arrow $w$. 
\end{notation}

In the study of Terracini loci (even not on curves) the paper \cite{CiMiranda98} is very useful and we explain it in the following remark used in the proof of Proposition \ref{ratiomax}.

\begin{remark}\label{oo1}
Fix an integral projective curve $X$ in $ \Pj^r$, and let $W$ be  the image of the restriction map $H^0(\Oc_{\Pj^r}(1)) \to H^0(\Oc_X(1))$.  For any zero-dimensional scheme $Z\subset X$ set $W(-Z):= H^0(\mathcal{I}_Z(1))\cap W$. Fix integers $s>0$ and $e_i>0$, $i=1,\dots ,s$. Let $Z$ be a general  subscheme of $X_{\reg}$ with $s$ connected components of length $e_1,\dots,e_s$.
Then we claim that the main result of \cite{CiMiranda98} and its proof imply that:
$$\dim W(-Z) =\max \{0,\dim W-e_1-\cdots -e_s\}.$$ 
Namely the proof there yields that a general scheme in $X$ (hence obviously curvilinear) which sits in no hypersurfaces of degree $d$  imposes independent conditions to the linear system of hypersurfaces of degree $d$. \\
Notice that the scheme $Z$  is a Cartier divisor of $X$, and the integer $\dim W(-Z)$ is the codimension of the linear space $\langle Z\rangle$ in $\Pj^r$.
\end{remark}

For rational curves we can easily show the following.

\begin{proposition}\label{ratio}
Fix integers $d \ge r\ge 3$ and $x$ with $r\leq 2x\leq d$. Then there is a smooth and non-degenerate rational curve 
$X\subset \Pj^r$ such that $\T(X,x) \ne \emptyset$.
\end{proposition}
\begin{proof}
Fix a hyperplane $H\subset \Pj^r$. Let $Z\subset H$ be a general union of  $x$ arrows and let $Z'$ be a general set of $d-2x$ points. 
By \cite[Theorem 1.6]{Perrin87} there is a smooth and non-degenerate rational curve $X\subset \Pj^r$ such that $(Z\cup Z') =X\cap H$. 
Let $S\subset H$ be the reduction of $Z$. Since $X$ is smooth and $2x\geq r$, then $S\in \T(X,x)$.
\end{proof}

Indeed, we have a characterization of rational normal curves in terms of Terracini loci.

\begin{proposition}\label{ratiomax}
Let $r\ge 3$ be an odd integer and $X\subset \Pj^r$ a smooth, connected and non-degenerate curve. Then $X$ is a rational normal
curve if and only if $\T(X,(r+1)/2)) =\emptyset$.
\end{proposition}
\begin{proof}
Set $x:=(r+1)/2$. \\
 The `` if '' part is true, because if $X$ is a rational normal curve each  zero-dimensional scheme $Z\subset X$ of degree $r+1$
is linearly independent.\\
Now assume $\T(X,x) =\emptyset$. Set $d:= \deg (X)$ and fix a general $S\subset X$ of cardinality $x-1$, say $S=\{p_1,\dots ,p_{x-1}\}$.  Let
$V$ be the linear span of the double scheme $(2S,X)$. For all positive integers $a_1,\dots ,a_{x-1}$ set $V(a_1,\dots ,a_{x-1}):= \langle a_1p_1+\cdots +a_{x-1}p_{x-1}\rangle$. Note that $V=V(2,\dots ,2)$.
Since $S$ is general in $X$, Remark \ref{oo1} applied to the curve $X$  gives $$\dim V(a_1,\dots ,a_{x-1}) = \min \{r,x-2 + \sum _{i=1}^{x-1} (a_ i-1)\}.$$ Thus $\dim V = r-2$. Since $\dim V =r-2$ and $\dim V(a_1,\dots ,a_{x-1}) =r-1$ if there is $j\in \{1,\dots ,x-1\}$
such that $a_j=3$ and $a_i=2$ for all $i\ne j$,
the scheme $V\cap X$ contains each $p_i$ with multiplicity $2$.
First assume $(V\cap X)_{\red} \ne S$ and take $o\in (V\cap X)_{\red} \setminus S$. Since $\dim \langle 2o\cup V\rangle \le \dim V+1$,  $S\cup \{o\}\in \T(X,x)$, a contradiction.
Thus  $(V\cap X)_{\red} =S$. Since we proved that each $p_i$ appear with multiplicy $2$ in the scheme-theoretic intersection $V\cap X$, we have $V\cap X =(2S,X)$, which has degree $2x-2=r-1$. Let $u$
denote the linear projection from $V$ to $\Pj^1$. Since $X$ is smooth, $u_{|X\setminus S}$ extends to a morphism $u': X\to \Pj^1$, and
 the degree of $u'$ is $d-r+1$. The assumption $\T(X,x) =\emptyset$ implies that $u'$ has no ramification
point, except possibly at the points of $S$. Fix $a\in S$, say $a=p_1$. The point $p_1$ is a ramification point of $u'$ only if $V(4,2,\dots ,2)$ is a hyperplane. This is false, because $\dim V(4,2,\dots ,2)=r$. Hence $u'$ has  no ramification points. This is possible only  if $d=r$, hence $X$ is a rational normal curve.\end{proof}

Since in the rest of the paper we will often argue by induction, taking the smoothing of nodal, reducible curves, we need some preliminary results
on normal bundles of reducible curves.

\begin{remark}
Let $X\subset \Pj^r$ be a reduced curve with only locally complete intersection singularities. 
$N_X$ is a vector bundle of rank $(r-1)$ on $X$ and 
$$\deg (N_X) =(r+1)\deg (X) +(r-1)(1-p_a(X)).$$ 
There is a map $\phi :T\Pj^r_{|X}\to N_X$ 
which is surjective outside $\mathrm{Sing}(X)$. Consider the  restriction to $X$ of the Euler's sequence of $T\Pj^r$:
\begin{equation}\label{eqr1}
0\to \Oc_X\to \Oc_X(1)^{\oplus (r+1)}\to T\Pj^r_{|X} \to 0
\end{equation}
Now assume $h^1(\Oc_X(1))=0$. Since $\dim X =1$, then $h^2(\Oc_X)=0$. Thus \eqref{eqr1} gives $h^1(T\Pj^r_{|X})=0$. \\
The map $\phi$ induces the following exact sequence of coherent sheaves on $X$:
\begin{equation}\label{eqr2}
0 \to \mathrm{Ker}(\phi) \to T\Pj^r_{|X}\to \mathrm{Im}(\phi)\to 0
\end{equation}
Since $\dim X =1$, $h^2(\mathrm{Ker}(\phi))=0$. Thus \eqref{eqr2} gives $h^1(\mathrm{Im}(\phi))=0$.\\
 If $X$ is smooth, then $N_X=\mathrm{Im}(\phi)$ and hence $h^1(N_X)=0$. \\
 Now assume $X$ singular. Consider the exact sequence
\begin{equation}\label{eqr3}
0 \to \mathrm{Im}(\phi)\to N_X \to N_X/\mathrm{Im}(\phi)\to 0
\end{equation}
Since $N_X/\mathrm{Im}(\phi)$ is supported by the finite set $\mathrm{Sing}(X)$, then $h^1(N_X/\mathrm{Im}(\phi))=0$. Thus \eqref{eqr3}
gives $h^1(N_X)=0$ even if the non-special curve is singular. \\
If $X$ is smooth and rational, then $h^1(\Oc_X)=0$. As above we obtain $h^1(T\Pj^r(-1)_{|X})=0$ and $h^1(N_X(-1)) =0$.
\end{remark}

In particular we study the case of curves in $\Pj^3$.

\begin{lemma}\label{P3lem}
Let $C\subset \Pj^3$ be a reduced curve and let $C'\subset\Pj^3$ be a smooth conic that meets $C$ at $i\leq 3$ points $p_1,\dots,p_i\in C_{reg}$ with $C\cup C'$ nodal at each $p_i$.
Let $H$ be the plane spanned by $C'$. If $i>1$, assume also that  $T_{p_1}C\nsubseteq H$, $T_{p_2}C\nsubseteq H$ and $T_{p_1}C\cap T_{p_2}C=\emptyset$. \\
1) If $i=1$ then $N_{C\cup C'|C'}$ is the direct sum of two  line bundles, one of degree $3$ and one of degree $4$.\\
2)  If $i=2$ then $N_{C\cup C'|C'}$ is the direct sum of two  line bundles, both of degree $4$.\\
3) If $i=3$, then $N_{C\cup C'|C'}$ is the direct sum of two  line bundles, one of degree $4$ and one of degree $5$.\\
In any case $h^1(N_{C\cup C'|C'}(-1))=0$.
\end{lemma}
\begin{proof}
Note that $N_{C'}\cong \Oc_{C'}(2)\oplus \Oc_{C'}(1)$, i.e. $N_{C'}$ is the direct sum of a degree $4$ line bundle and a degree $2$ line bundle. \\
We prove 2) first. Set $Y:= T_{p_1}C\cup T_{p_2}C\cup C'$.
By \cite[Corollary 3.2]{HartHir85} the bundle $N_{C\cup C'|C'}$ is obtained from $N_{C'}$ by making two positive elementary
transformations  (in the sense of \cite[\S 2]{HartHir85}), which depend uniquely on  $C',p_1,p_2$, and the lines $T_{p_1}C,T_{p_2}C$. 
Thus $N_{C\cup C'|C'}\simeq N_{Y|C'}$. Since $Y$ is connected and $\deg(Y)=4$, $h^0(\Oc_Y(2))<10$. Thus $h^0(\Ic_Y(2))>0$.
Since  $T_{p_1}C\cap T_{p_2}C=\emptyset$, $Y$ is contained neither in a reducible quadric, nor in a quadric cone. Thus $Y$ lies in a smooth quadric $Q$ with, say, $T_{p_1}C\in |\Oc_Q(1,0)|$. Then also $T_{p_2}C$, which does not meet $T_{p_1}C$, lies in $ |\Oc_Q(1,0)|$. Obviously $C'\in |\Oc_Q(1,1)|$, and we conclude that  $Y\in |\Oc_Q(3,1)|$. Thus $N_Y$ is an extension of $\Oc_Y(2)$ by the restriction to $Y$ of the line bundle $\Oc_Q(3,1)$. Hence $N_{Y|C'}$ is an extension of two line bundles of degree $4$.\\
1) By \cite[Corollary 3.2]{HartHir85} the bundle $N_{C\cup C'|C'}$ is obtained from $N_{C'}$ by making one general positive elementary
transformation, and as above one general positive elementary transformation of $N_{C\cup C'|C'}$ 
 is a direct sum of a line bundle of degree $4$ and a line bundle of degree $3$. \\
3) Here $N_{C\cup C'|C'}$ is obtained from $N_{C'}$ by making $3$ positive elementary
transformations, one at $p_1$, one at $p_2$, and one at $p_3$. Hence $N_{C\cup C'|C'}$ is obtained from one positive elementary transformation
on   the vector bundle $E$ on $C'$ obtained by making positive elementary transformations only at $p_1$ and $p_2$. By 1) $E$ is isomorphic 
to $\Oc_{C'}(2)\oplus\Oc_{C'}(2)$. Since any vector bundle on $C'$ splits in  a direct sum of line bundles, the claim follows.\end{proof}

\begin{remark}\label{Kleppe} We will use repeatedly a result of J. Kleppe, who
  proved  that the functor of deformations inside $\Pj^r$ of a (possibly reducible) curve $T\subset \Pj^r$ with fixed  hyperplane section 
$H\cap T$ has tangent space $H^0(T,N_T(-1))$ and obstruction space $H^1(N_T(-1))$. See e.g. \cite[Th. 1.5]{Perrin87}.
\end{remark}

\begin{remark}\label{aggiunto1}
Let $C$ and $D$ be locally complete intersection space curves. Assume $S:= C\cap D=C_{\reg}\cap D_{\reg}$ and that at each $p\in S$ the tangent lines of $C$ and $D$ at $p$ are different. Set $Y:= C\cup D$.
Call $N_C^+$ (resp. $N_D^+$) the rank $2$ vector bundle on $C$ (resp. $D$) obtained from $N_C$ (resp. $N_D$) making a positive elementary transformation in the direction of $T_pD$ (resp. $T_pC$)
at all $p\in S$. We have $\deg (N_C^+)=\deg (N_C)+\#S$ and $\deg (N_D^+) =\deg (N_D)+\#S$. By \cite[Cor. 3.2]{HartHir85} we have $N_C^+ \cong N_{Y|D}$ and $N_D^+ \cong N_{Y|D}$. Thus the Mayer-Vietoris exact sequence of the vector bundle $N_Y$ gives  the following  exact sequence
\begin{equation}\label{eqagg1}
0 \to N_{C\cup D} \to N_C^+\oplus N_D^+ \to N_{Y|S}\to 0.
\end{equation}
Thus to prove that $h^1(N_Y(-1)) =0$  it is sufficient to prove that $h^1(N_C^+(-1)) = h^1(N_D^+(-1)) =0$
and that the restriction map $H^0(D,N_D^+(-1)) \to H^0(N_{Y|S})$ is surjective.

A similar statement holds for the proof that $h^1(N_{C\cup D})=0$.
\end{remark}

\begin{lemma}\label{sab3}
Let $H\subset \Pj^3$ be a plane and let $C'\subset H$ be a smooth conic. Fix $3$ distinct points $p_1,p_2,p_3$ of $C'$ and lines $L_1,L_2,L_3$ such that $H\cap L_i =\{p_i\}$ for all $i$.
Let $E$ (resp. $F$) be the vector bundle obtained from $N_{C'}$ making  positive elementary transformations at $p_1$ and $p_2$ (resp. $p_1$, $p_2$ and $p_3$) with respect to $L_1$ and $L_2$ (resp. $L_1$, $L_2$ and $L_3$). Then $E$ is a direct sum of $2$ line bundles of degree $4$ and $F$ is a direct sum of a line bundle of degree $5$ and a line bundle of degree $4$.
\end{lemma}

\begin{proof}
Set $X:= C'\cup L_1\cup L_2$ and $Y:= X\cup L_3$. Since $H\cap L_i=\{p_i\}$, then $L_i$ is transversal to $H$. Thus $X$ and $Y$ are nodal at $p_1$, $p_2$ and $p_3$. Note that $E\cong N_{X|C'}$ and $F\cong N_{Y|C'}$. 

First assume $L_1\cap L_2\ne \emptyset$. Thus $X$ is nodal with $3$ nodes and arithmetic genus $1$. Call $M$ the plane containing $L_1\cup L_2$. To prove that $E\cong \mathcal{O}_{C'}(2)\oplus \mathcal{O}_{C'}(2)$ it is sufficient to prove
that $X$ is the complete intersection of $M\cup H$ and a quadric. This is true, because  $h^0(\mathcal{O} _X(2)) =8$, $h^0(\mathcal{O}_{\Pj^3}(2)) =10$ and hence $h^0(\mathcal{I} _X(2)) \ge 2$.

Now assume $L_1\cap L_2=\emptyset$. In this case $X$ is contained neither in a reducible quadric nor in a quadric cone. Since $h^0(\mathcal{O}_X(2)) =9$, $X$ is contained in a smooth quadric, $Q$.
Call $|\mathcal{O}_Q(1,0)|$ the ruling of $Q$ containing $L_1$, and hence also containing $L_2$. Since $C'\in |\mathcal{O}_Q(1,1)|$, then $X\in |\mathcal{O}_Q(3,1)|$. Since $N_Q\cong \mathcal{O}_Q(2)$, we have an exact sequence
\begin{equation}\label{eqsab2}
0 \to N_{X,Q} \to N_X \to N_X(2)\to 0
\end{equation}
We have $N_{X.Q} \cong \mathcal{O}_X(3,1)$ so that its restriction to $C'$ has degree $4$. Since \eqref{eqsab2} is an exact sequence of vector bundles, its restriction to $C'$ is an exact sequence of vector bundles on $C'\cong \Pj^1$ in which the leftmost and the rightmost terms are
line bundles of degree $4$. Since $C'\cong \Pj^1$, $E$ is a direct sum of two line bundles of degree $4$.

The bundle $F$ is obtained from $E$ making a positive elementary transformation, and all rank $2$ vector bundles on $\Pj^1$ split. Hence $F$ is a direct sum of a line bundle of degree $5$ and a line bundle of degree $4$.\end{proof}

Now we are ready to prove some non-emptiness results in $\Pj^3$. \\
Recall that  $H(d,g,3)$ denotes the set of smooth and non-degenerate curves
space curves $X$ of degree $d$ and genus $g$, such that $h^1(\Oc_X(1))=0$.

\begin{theorem}\label{i1}
Fix integers $g\ge 0$, $d\ge g+3$ and $x$ such that $6\le 2x\le d$. Then there is $X\in H(d,g,3)$ such that $\T(X,x)\ne\emptyset$
and $h^1(N_X(-1)) = 0$.
\end{theorem}
\begin{proof}
Fix a plane $H\subset \Pj^3$. We will find $X\in H(d,g,3)$ such that $X$ is tangent to $H$ at $x$ points of $H$ spanning $H$.
We first dispose of the case $x =\lfloor d/2\rfloor$ and $d=g+3$ in steps (a) and (b), leaving the case
$x=\lfloor d/2\rfloor$ and $d>g+3$ to step (c) and the case $3\le x < \lfloor d/2\rfloor$ to step (d).

 (a) Assume $d$ even and $x=d/2$. 
 
\quad (a1) Assume $d=6$. By Proposition \ref{ratiomax} $\T(C,2)\ne \emptyset$ for each smooth curve $C$ of genus $1$ and degree $4$.
Hence by a change of coordinates we find a smooth curve $C\subset \Pj^3$ of degree $4$ and genus
$1$ which is tangent to $H$ at $2$ distinct points, say $q_1$ and $q_2$. Since $C$ is the complete intersection of $2$
quadric surfaces, the normal bundle of $C$ splits as $N_C\cong \Oc_C(2)\oplus \Oc_C(2)$. 
Since $C$ has genus $1$, we get $h^1(N_C(-1))=0$. 
Fix a general $q_3\in H$ and let $M\subset \Pj^3$ be a general plane containing $q_3$. Since $q_3$ and $M$ are
general, we may assume that $M$ is transversal to $C$, $C\cap M \cap H =\emptyset$, and $C\cap M$ spans $M$. 
Fix three distinct points $\{p_1,p_2,p_3\}\in C\cap M$ which span $M$. There is a smooth conic $C'$ containing
$\{p_1,p_2,p_3,q_3\}$, tangent in $q_3$ to the line $H\cap M$  but not containing the fourth point $p_4$ of $C\cap M$. 
Indeed, $q_3$ is general in $H\cap M$, while there are only at most two conics passing through $p_1, p_2, p_3, p_4$ 
and tangent to $H\cap M$. Since $M$ is transversal to $C$,  the curve $Y:= C\cup C'$ 
is nodal. By Lemma \ref{P3lem} the vector bundle $N_{C'}^+(-1)$ is the direct sum of a line bundle of degree $2$ and a line bundle of degree $3$.
Since line bundles of degree $2,3$ separate any set of three points in the smooth conic $C'$, the restriction map $H^0(C',N_{C'}^+(-1)) \to H^0(N_{C'}^+(-1)_{|\{p_1,p_2,p_3\}})$ is surjective. 
Remark \ref{aggiunto1} gives $h^1(N_Y(-1))=0$, so that $Y$ is smoothable by \cite[Theorem 4.1]{HartHir85}. 
By semicontinuity,  a general member $X_0$ of a smoothing family 
of $Y$ satisfies $h^1(N_{X_0} (-1)) = 0$.
By construction $\T(Y,3)\ne \emptyset$.  Yet, in order to conclude, we need more: we need to smooth $Y$ in a family of space curves 
whose elements $Y_\lambda$ satisfies $\T(Y_\lambda,3)\ne \emptyset$. In other words, we need:

\quad {\bf Claim 1:} There are an affine smooth and connected curve $\Delta$, $o\in \Delta$ and a flat family 
$\{Y_\lambda\}_{\lambda\in \Delta}$ of space curves such that $Y_0=Y$, the general element of  the family is smooth, 
and $\T(Y_\lambda,3)\ne \emptyset$ for all $\lambda\in \Delta$.

\quad {\bf Proof of Claim 1:} Set $Z':= (2q_1,C)\cup (2q_2,C)$, $Z'':= (2q_3,C')$ and $Z:= Z'\cup Z''$. 
Note that $Z\cap C =Z'$ and $Z\cap C'=Z''$. Since $q_1$, $q_2$ and $q_3$ are smooth points of $Y$, $Z$ is a degree $6$ Cartier divisor of $Y$. Thus $N_Y(-Z)$ is a rank $2$ vector bundle on $Y$ with $\deg (N_Y(-Z)) =\deg (N_Y)-12$.
The vector space $H^0(N_Y(-Z))$ is the tangent space to the functor 
of deformations of $Y$ inside $\Pj^3$ in families of curves containing $Z$, while $H^1(N_Y(-Z))$ is an obstruction space of this 
functor \cite[Th. 1.5]{Perrin87}. To prove Claim 1 it is sufficient to find a smoothing family $\{Y_\lambda\}_{\lambda\in \Delta}$ of $Y=Y_o$
such that $Z\subset Y_\lambda$ for all $\lambda \in \Delta$. Note that $N_{C'}(-Z) =N_{C'}(-Z'')$ is a direct 
sum of degree $2$ line bundles  on $C'\cong \Pj^1$. Thus $h^1(N_{C'}(-Z-p_i)) =0$ for all $i$. 
In particular (with the terminology of \cite{HartHir85})
$h^1(N_C(-Z)^-) =0$. Since $N_C \cong \Oc_C(2)\oplus \Oc_C(2)$, and $\deg (Z\cap C)=4$, $h^1(N_C(-Z-p_1-p_2-p_3)) =0$. 
Thus $h^1(F)=0$ for every vector bundle $F$ on $C$ obtained from $N_C(-Z-p_1-p_2-p_3)$ making finitely many 
positive elementary transformations. Consider the Mayer-Vietoris exact sequence
\begin{equation}\label{eqbo1}
0 \to N_Y(-Z) \to N_Y(-Z)_{|C} \oplus N_Y(-Z)_{|C'}\to N_Y(-Z)_{|\{p_1,p_2,p_3\}}\to 0
\end{equation}
Since $N_Y(-Z)_{|C}(-p_1-p_2-p_3)$ is obtained from $N_C(-Z-p_1-p_2-p_3)$ by making positive elementary transformations,   
$h^1(N_Y(-Z)_{|C}(-p_1-p_2-p_3)) =0$
(here we consider $p_1,p_2,p_3$ as points of the smooth curve $C$). The sequence
$$ 
0 \to N_Y(-Z)_{|C} (-p_1-p_2-p_3) \to N_Y(-Z)_{|C} \to N_Y(-Z)_{|\{p_1,p_2,p_3\}}\to 0
$$
 shows that  the  map $\phi:H^0(N_Y(-Z)_{|C})\to H^0(N_Y(-Z)_{|\{p_1,p_2,p_3\}})$ is surjective, and $h^1(N_Y(-Z)_{|C})=0$. We also know that $N_Y(-Z)_{|C'}$ is a direct sum of $2$ line bundles of non-negative degree on the smooth conic $C'$, so that $h^1(N_Y(-Z)_{|C'})=0$. Hence the surjectivity of $\phi$ gives a fortiori, in sequence \eqref{eqbo1}, that  $h^1(N_Y(-Z))=0$. Thus we may apply the proof of \cite[Th. 4.1]{HartHir85} since the deformation functor of $Y$ which maintains $Z$ fixed is unobstructed. We obtain a family  $\{Y_\lambda\}_{\lambda\in \Delta}$ as in the statement, in which the general element contains $Z$, hence it is tangent to $H$ at three points.
\hfill\qed

\quad (a2) Assume $d\ge 8$ and that the theorem is true for the triples $(d',g',x')$ such that $x=d'/2$, $d'=g'+3$ and $d'\le d-2$. 
Take a solution $C$ for $(d',g',x') =(d-2,g-2,x-1)$ and use the proof of step (a1) with this $C$ instead of an elliptic curve,
as follows. By induction, there is a plane $H$ which contains a set $Z'$ of $x-1$ arrows in $C$. Then take a general plane 
$M$ and three points  $o_1,o_2,o_3$ of $M\cap C$. Fix a conic $C'\subset M$ passing through $o_1,o_2,o_3$ and tangent to $H$ at a point $p$. Since the base locus of the family of conics in $M$ passing through $o_1,o_2,o_3$ and tangent to $H$ is exactly given by $o_1,o_2,o_3$, we may assume that $C'$ misses any other point of $M\cap C$. Set $Y=C\cup C'$, $Z''=(2p,C')$, and $Z=Z'\cup Z''$.

\quad {\bf Claim 2:} $h^1(N_Y(-Z)) =0$. 

\quad {\bf Proof of Claim 2:} The bundle $N_{C'}^+(-Z''-o_1-o_2-o_3)$ is a direct sum of a line bundle of degree $0$ and a line bundle of degree $-1$. Thus $h^1(C',N_{Y|C'}(-Z-o_1-o_2-o_3))=0$.
As in Claim 1, just working on $C'$ instead of $C$ and with $o_1,o_2,o_3$ instead of $p_1,p_2,p_3$, this implies  that the restriction map  $$H^0(N_Y(-Z)_{|C'})\to H^0(N_Y(-Z)_{|\{o_1,o_2,o_3\}})$$ is surjective. Then the analogue of  \eqref{eqbo1} gives $h^1(N_Y(-Z)) =0$.
 \hfill\qed
 
So, we can continue the induction by taking a general element $X$ in a family $\{Y_\lambda\}_{\lambda\in \Delta}$ giving a smoothing  of $Y$
 and fixing $Z$. Notice that at any step $X$  satisfies $h^1(N_X(-Z)) =0$ by semicontinuity, since $h^1(N_Y(-Z)) =0$.

\quad (b) Assume $d$ odd and $x = (d-1)/2$. Since $x\geq 3$, $d\ge 7$. By Proposition \ref{ratiomax} each $C\in H(5,2,3)$ has
$\T(C,2)\ne \emptyset$. To adapt the proof of step (a) we first find $C\in H(5,2,3)$ with $h^1(N_C(-1))=0$. Let
$Y'\subset \Pj^3$ be a smooth rational cubic. Take a plane $M\subset \Pj^3$ transversal to $Y'$ and let $C'\subset M$ 
be a smooth conic containing $Y'\cap M$. Set $T:= Y'\cup C'$. Lemma \ref{P3lem}
 gives $h^1(N_T(-1))=0$. By \cite[Theorem 4.1]{HartHir85} and Lemma \ref{P3lem}, $T$ is smoothable in a family preserving  the scheme-theoretic intersection $(Y'\cup C')\cap H$ as in step (a1). By semicontinuity, a general element $C$ of the family has $h^1(N_C(-1))=0$.
 Next we proceed by induction as in step (a), adding suitable conics $C'$, constructed as in step (a), to curves $C$ which are solutions for the triple $(d',g',x') =(d-2,g-2,x-1)$, and taking a general smoothing $X$ of $C\cup C'$ . To see that we get $h^1(N_X(-1)) =0$ we apply Lemma \ref{sab3} to $C'$.

\quad (c) Assume $d>g+3$.  If $d-g-3$ is even, start with a curve of genus $g$ and degree $d'=g+3$, constructed as in step (a2) when $d'$ is even, or constructed as in
step (b) if $d'$ is odd. Continue for $ (d-g-3)/2$ steps, by adding to the previously constructed curve $C$ a smooth conic $C'$,  tangent to $H$, with $\#(C'\cap C)=1$ and $C'\cup C$ nodal.  By Lemma \ref{P3lem} we always get $h^1(N_{C\cup C'}(-1))=0$. After  $ (d-g-3)/2$ steps we get the claimed curve.  \\
Assume that $d-g-3$ is odd. When $g=0$, start  with a rational quartic, by using Proposition \ref{ratiomax}.  When $g\geq 1$, start with a curve $C$ of genus $g-1$ and degree $d'=g+2$, constructed as in step (a2) when $d'$ is even, or constructed as in
step (b) if $d'$ is odd. In the first step  add  to $C$ a smooth conic $C'$, tangent to $H$,
which meets the curve $C$ at two points whose tangent lines $t_1,t_2$ to $C$  are disjoint and different from the tangent lines to $C'$. We may take such a conic $C'$ in 
 a plane which does not contain $t_1,t_2$. Then  continue for $ (d-g-4)/2$ steps, by adding to the previously constructed curve $C$ a smooth conic $C'$, tangent to $H$, which intersects $C$ at a unique point and with $C\cup C'$ nodal. In any case, the assumptions of Lemma \ref{P3lem} hold for $C\cup C'$. Thus by Lemma \ref{P3lem} we always get $h^1(N_{C\cup C'}(-1))=0$, and we can continue the induction, by taking a smoothing of $C\cup C'$ which preserves the intersection with $H$.

\quad (d) Assume $3\le x \le \lfloor d/2\rfloor$. We start with some curve $Y$   such that $\T(Y,2)\ne \emptyset$ and $h^1(N_Y(-1))=0$. We take $Y$ of genus $1$ and degree $4$ if $d$ is even or genus $2$ and degree $5$ if $g$ is odd. Then we continue as in steps (a), (b) and (c) above,   except that in $\lfloor d/2\rfloor-x$ steps  we add a smooth conic $C'$ not tangent to $H$.
\end{proof}

For space curves $X\subset \Pj^3$ with $h^1(\Oc_X(1)) \ne 0$ we prove the following result.

\begin{theorem}\label{i1.1}
Fix integers $x\ge 3$, $g\ge 0$ and $d$ such that $g\ge 2x+2$. Write $g=2x+1+5s-q$, $0\le q\le 4$, and assume $d\ge 2x+4+3s$. Then there is a smooth and connected curve
$X\subset \Pj^3$ such that $\deg (X) =d$, $p_a(X)=g$, $h^1(N_X)=0$, $h^1(\Oc_X(2)) =0$ and $\T(X,x)\ne \emptyset$. 
\end{theorem}

\begin{remark}\label{ll0}In Theorem \ref{i1.1}, for each fixed $x$ we find $g_0$ such that for all $g\ge g_0$ and all $d\ge \frac{3}{4}g+3$ 
 there is a smooth curve $X\subset \Pj^3$ of genus $g$ and degree $d$ with $\T(X,x)\ne \emptyset$.
The same is true for a slowing increasing function $x(g)$ of $g$.  \\
Note that for a fixed $x$ and for $g\gg x$ these curves $X$ cover a range of degrees and genera larger that the Brill-Noether
range $d\ge \frac{3}{4}g+3$.
\end{remark}

For the proof of Theorem \ref{i1.1}, we need a series of preliminary lemmas.

\begin{lemma}\label{3a2}
Fix an integer $e\in \{1,2,3,4\}$. Let $C\subset \Pj^3$ be an integral and non-degenerate curve of degree $d$. If $e=4$ assume $d\ge 4$. Take a union $Y\subset \Pj^3$ of finitely many curves such that $C\nsubseteq Y$. Then there is a smooth conic $C'$
such that $\#(C\cap C')=e$, $C\cup C'$ is nodal at each point of $C\cap C'$, and $C'\cap Y=\emptyset$.
\end{lemma}

\begin{proof}
Take a general plane $M\subset \Pj^3$. The plane $M$ is transversal to $C$, $Y\cap M$ is finite and $Y\cap C\cap M=\emptyset$. By the trisecant lemma no $3$ of the $d$ points of $C\cap M$ are collinear.
Fix $S\subseteq C\cap M$ such $\#S=e$. Since no $3$ of the points of $S$ are collinear, $S$ is the scheme-theoretic base locus of the $(5-e)$-dimensional linear space  $|\mathcal{I}_S(2)|$. A general element of $|\mathcal{I}_S(2)|$
is smooth. Thus there is a smooth $C'\in |\mathcal {I}_S(2)|$ such that $C'\cap Y=\emptyset$ and $C'\cap C=S$. Since $M$ is transversal to $C$ and $C'\subset M$, then $C\cup C'$ is nodal at each point of $S$.
\end{proof}
\smallskip

\begin{lemma}\label{ai6}
Let $C\subset \Pj^3$ be an integral and non-degenerate curve of degree $4$. For a general $S\subset C$ of length $6$ there is a
rational normal curve $T_S\subset \Pj^3$ such that $S=C\cap T_S$ and $T_S\cup C$ is nodal.\\
 Moreover, if $p_a(C)=1$ then $T_{S_1}\cap T_{S_2}=\emptyset$ for a general $S_1\times S_2\subset C\times C$ such that both $S_1$ and $S_2$ have length $6$.
\end{lemma}

\begin{proof}
Let $\mathcal {U}$ denote the set of all $A\subset \Pj^3$ such that $\#A=6$ and $A$ is in linear general position. For every $A\in \mathcal {U}$ there is a unique rational normal curve $T_A$ containing $A$.
Set $\mathcal {U}(C):= \{S\in \mathcal {U}\mid S\subset C\}$. Since $C$ is integral and non-degenerate, $\mathcal{U}(C)$ is an integral quasi-projective variety of dimension $6$. The set
$S$ is a general element of $\mathcal{U}(C)$.  Since $\deg (C)\ne 3$, then $C\ne
T_S$. We need to prove that $S$ is equal to $C\cap T_S$ (set-theoretically), that $T\cup C$ is nodal, and the last
assertion, concerning a general $S_1\times S_2\subset C\times C$, when $p_a(C)=1$.  \\
If $p_a(C)=1$, the curve $C$ is the complete intersection of $2$ quadrics and (since it has at most one singular point and with embedding
dimension $2$) it is contained in a smooth quadric $Q$, say $C\in |\Oc_Q(2,2)|$. Fix a general $S\subset C$ of length $6$. 
By the generality  of $S$, since $h^0(\Oc_Q(1,2))=h^0(\Oc_Q(2,1))=6$, we get $h^0(\Ic _{S,Q}(2,1)) = h^0(\Ic_{S,Q}(1,2))=0$. Thus $T_S\nsubseteq Q$. Bezout theorem gives $S =T_S\cap Q$ as schemes. Thus $T_S\cap C =S$ and $T_S\cup C$ is nodal. 
If $p_a(C)=0$ then $h^0(\Ic_C(2))=1$ and the unique quadric surface $Q$ containing
$C$ is smooth \cite[Ex. V.2.9]{Hartshorne}, with either $C\in |\Oc_Q(1,3)|$ or $C\in |\Oc_Q(3,1)|$. We conclude as in the case $p_a(C)=1$.\\
Now we prove the last claim, for $p_a(C)=1$. It is sufficient to find $S_1,S_2\in \mathcal{U}(C)$ such that $T_{S_1}\cap T_{S_2} =\emptyset$. 
The pencil $|\mathcal{I} _C(2)|$ has only finitely many (i.e. $4$) singular elements. 
Take smooth quadrics $Q_1,Q_2\in |\mathcal{I}_S(2)|$ such that $Q_1\ne Q_2$.
Note that $C\in |\mathcal{O}_{Q_i}(2,2)|$, $i=1,2$. Take a general $T_i\in |\mathcal{O}_{Q_i}(2,1)|$. Each $T_i$ is a rational normal curve and $\deg (T_i\cap C) =6$. Bertini Theorem gives that $S_i:= T_i\cap C$ has cardinality $6$.
Since $S_i\subset T_i$ and $T_i$ is a rational normal curve, $S_i$ is in linear general position. Thus $T_i =T_{S_i}$. Since $T_{S_i}\ne C$, we have $T_{S_1}\nsubseteq Q_2$ and $T_{S_2}\nsubseteq Q_1$.
Since $T_i\subset Q_i$, $i=1,2$, we get $T_{S_1}\cap Q_2 =T_{S_1}\cap Q_1\cap Q_2 =T_{S_1}\cap C =S_1$. Since $T_{S_2}\subset Q_2$, then $T_{S_1}\cap T_{S_2}\subseteq T_{S_1}\cap Q_2 =S_1$. Since $S_1\cap S_2=\emptyset$ and $T_{S_2}\cap C =S_2$, we get $T_{S_1}\cap T_{S_2}=\emptyset$.
\end{proof}

\begin{lemma} \label{x0} For $q=0,\dots,4$ and $x\geq 2$ the space $H(2x+4, 2x+1-q,3)$ contains smooth curves $X_0$ which satisfy the following conditions.
\begin{enumerate}
\item $\T(X_0,x)$ is non-empty, and there exists a plane $H$ which contains a union of $x$ arrows $Z\subset X_0$ with $h^1(N_{X_0}(-Z))=0$.
\item For a general choice of $s$ subsets $S_1,\dots, S_s\subset X_0$ of cardinality $6$, there are disjoint rational normal curves $T_1,\dots, T_s$ such that for all $i$ the union $X_0\cup T_i$ is nodal, and $T_i$ meets $X_0$ exactly at $S_i$. 
\end{enumerate}
\end{lemma}
\begin{proof}
We make induction on $x\geq 2$.\\
To deal with the case  $x=2$ (i.e. $d=8$ and $g=5-q$) we start with a smooth elliptic quartic curve $C$. We know that $H^1(N_C(-1))=0$ and $\T(C,2)$ is non-empty. 
Since $\T(C,2)\neq \emptyset $, there is a subset $Z\subset C$ of cardinality $2$ such that the tangent lines to $C$ at the points of $Z$ lie in a plane $H$. By Lemma \ref{ai6}, for a general choice of $s$ subsets $S_1,\dots, S_s\subset C$ of cardinality $6$ we find disjoint rational normal curves $T_1,\dots, T_s$ such that $C\cup T_i$ is nodal and $T_i\cap C=S_i$ for all $i$. Take a general plane $H'\subset \Pj^3$. Thus $H'$ is transversal to $C$  and $Z\cap H'=\emptyset$.
Fix $e\in \{0,1,2,3\}$ and $W\subset H'\cap C$ such that $\#W =4-e$.
Take a general smooth conic $C'\subset H'$ containing $W$. Since $\#W\le 4$ and $C'$ is general, $T_i\cap C'=\emptyset$ for all $i$.
Set $Y:= C\cup C'$. The normal bundle $N_Y(-Z)_{|C}$ is obtained from $N_C(-Z)$ by making positive elementary transformations, so that
$h^1(N_Y(-Z)_{|C}) =0$. Recall that $1\le \#W\le 4$. For each $o\in W$ the tangent line $T_oC$ of $C$ at $o$ is transversal to $H'$. If $\#W=4$ Lemma \ref{sab3} shows that the restriction map $\rho: H^0(N_Y(-Z)_{|C'})\to H^0(N_Y(-Z)_{|W})$ surjects, so that the exact sequence
\begin{equation}\label{eqbo2}
0 \to N_Y(-Z) \to N_Y(-Z)_{|C} \oplus N_Y(-Z)_{|C'}\to N_Y(-Z)_{|W}\to 0
\end{equation}
proves that $h^1(N_Y(-Z))=0$. If $\#W\le 3$, then the surjectivity of $\rho$ only use that $N_{C'}$ is a direct sum of $2$ line bundles of degree $\ge 2$ and that $N_Y(-Z)_{|C'}$ is obtained from $N_{C'}$ making positive elementary transformations.
By semicontinuity, we can find a smoothing $D$ of $C\cup C'$  which preserves $Z$. Thus we get  curves $D$ in $H(6,g,3)$, for $g = \#W \in\{1,2,3,4\}$, such that $\T(D,2)\neq \emptyset$, $h^1(N_D(-Z))=0$, and such that for general subsets $S_1,\dots, S_s\subset X_0$ of cardinality $6$ we have the disjoint rational normal curves $T_1,\dots, T_s$, as in the statement. 
Then consider a conic $C''$ which meets $D$ in one or in two general points. As above, in both cases $Y'=C''\cup D$ satisfies $h^1(N_{Y'}(-Z))=0$, and we can find a general smoothing $X_0$ of $Y'$ which preserves $Z$. $X_0$ has degree $8$, and for its genus $g$ we can obtain any number between $1$ and $5$, because $\#W$ is any integer between $1$ and $4$. Moreover, $h^1(N_{X_0}(-Z))=0$, $\T(X_0,2)\neq \emptyset$, and by semicontinuity $X_0$ satisfies condition 2 of the statement. This concludes the case $x=2$.\\
Assume we constructed the required curve $X\in H(2(x-1)+4, 2(x-1)+1-q,3)$. There exists a subscheme $Z\subset X$ formed by $x-1$ arrows which is contained in a plane $H$, and moreover $h^1(N_{X}(-Z))=0$. Take a general plane $M$ which is transversal to both $H$ and $X$, and misses $Z$. As in the proof of (a1) of  Theorem \ref{i1}, for $p_1,p_2,p_3\in M\cap X$ there exists a smooth conic $C'$ passing through $p_1,p_2,p_3$, tangent to $H$ at $p_4\notin X$, which misses the remaining points of $M\cap X$. Take $Y=X\cup C'$, and define $Z_0$ as the union of $Z$ and the arrow at $p_4$ tangent to $C'$. Notice that $Z_0$ lies in $H$. As in the proof of Claim 1, the vanishing of $H^1(N_X(-Z))$ implies the surjectivity of the map $H^0(N_Y(-Z_0)_{|X})\to H^0(N_Y(-Z_0)_{|\{p_1,p_2,p_3\}})$, and the analogue of sequence \eqref{eqbo1} shows that $h^1(N_{Y}(-Z_0))=0$. Thus there exists a smoothing $X_0\in H(2x+4,2x+1-q ,3)$ of $Y$ which preserves $Z_0$. The existence of $Z_0\subset X_0$ implies $\T(X_0,x)\neq \emptyset$. By semicontinuity $h^1(N_{X_0}(-Z_0))=0$, and $X_0$ satisfies condition 2 of the statement.
\end{proof}

We notice that, in the previous lemma, the condition  $h^1(N_{X_0}(-Z))=0$ implies $h^1(N_{X_0})=0$.
\medskip

\begin{lemma}\label{eae1}
Take a reduced curve $C\subset \Pj^3$ such that $h^1(\mathcal{O}_C(2))=0$. Let $T\subset \Pj^3$ be a smooth rational curve which is not an irreducible component of $C$. Set $f:= \deg(T)$
and $e:= \deg (C\cap T)$. If $e\le 2f+1$, then $h^1(\mathcal{O}_{C\cup T}(2)) =0$.
\end{lemma}

\begin{proof}
Consider the Mayer-Vietoris exact sequence
\begin{equation}\label{eqeae1}
0\to \mathcal{O}_{C\cup T}(2)\to \mathcal{O}_C(2)\oplus \mathcal{O}_T(2)\to \mathcal{O}_{T\cap C}(2)\to 0
\end{equation}
By assumption $h^1(\mathcal{O}_C(2))=0$. Obviously, $h^1(\mathcal{O}_T(2))=0$. Since $e\le 2f+1$, $T\cong \Pj^1$ and $\deg (\mathcal{O}_T(2))=2f$, the restriction map $H^0(\mathcal{O}_T(2))\to H^0(\mathcal{O}_{C\cap T}(2))$ is surjective.
Hence \eqref{eqeae1} gives $h^1(\mathcal{O}_{C\cup T}(2)) =0$.
\end{proof}

\begin{lemma}\label{ooo2}
Let $C\subset \Pj^3$ be an integral locally complete intersection curve and $Z\subset C_{\reg}$ a zero-dimensional scheme. Assume $h^1(N_C(-Z))=0$.
Take $S\subset C_{\reg}$ such that $\#S=6$, $S\cap Z=\emptyset$, $S$ is in linearly general position and the only rational normal curve $T$ containing $S$ meets $C$ only at $S$ and with $C\cup T$ nodal at each point of $S$. Set
$Y:= C\cup T$. Then $h^1(N_Y(-Z)) =0$.\end{lemma}

\begin{proof}
By assumption $S$ is the scheme-theoretic intersection of $C$ and $T$. Thus we have the following Mayer-Vietoris exact sequence
\begin{equation}\label{eqooo1}
0 \to N_Y(-Z)\to N_Y(-Z)_{|C}\oplus N_Y(-Z)_{|T}\to N_Y(-Z)_{|S}\to 0
\end{equation}
Since $Z\cap T=\emptyset$, $N_Y(-Z)_{|T}\cong N_{Y|T}$. Since $N_Y(-Z)_{|C}$ is obtained from $N_C(-Z)$ making positive elementary transformations and $h^1(N_C(-Z))=0$, $h^1(N_Y(-Z)_{|C})=0$. Since $N_T$ is a direct sum of $2$ line bundles of degree $5$, $\#S=6$ and
$N_Y(-Z)_{|T}$ is obtained from $N_T$ making positive elementary transformations, $h^1(N_Y(-Z)_{|T})=0$ and the restriction map $H^0(N_Y(-Z)_{|T})\to N_Y(-Z)_{|S}$ is surjective. Thus \eqref{eqooo1} gives $h^1(N_Y(-Z)) =0$.
\end{proof}

Now we are ready for the proof of Theorem \ref{i1.1}.
\medskip

{\it Proof of Theorem \ref{i1.1}.} 
Recall that $g=2x+1+5s-q$, where $q\in \{0,\dots,4\}$. By assumption $d\geq 2x+4+3s$.\\
We first dispose of the case $d= 2x+4+3s$. The case $s=0$ is covered by Lemma \ref{x0}. 
Consider a curve $X_0\in H(4x+2,2x+1,3)$ as in the statement of  Lemma \ref{x0}. Since $\T(X_0,x)$ is non-empty, we can take a set $Z$ of $x$ coplanar arrows in $X_0$ supported at $x$ points, with $h^1(N_{X_0}(-Z))=0$. Moreover, for a choice of $s$ general subsets $S_1,\dots, S_s\subset X_0$ of cardinality $6$ there are disjoint rational normal curves $T_1,\dots, T_s$ such that for all $i$ the union $X_0\cup T_i$ is nodal, and $T_i$ meets $X_0$ exactly at $S_i$.  Define $Y=X_0\cup T_1\dots\cup T_s$. Then $h^1(N_Y(-Z))=0$. Arguing by induction on $s$  one finds a smoothing $X_s$ of $Y$  which preserves $Z$. The curve $X_s$ belongs to $H(2x+4+3s, 1+2x+5s-q,3)$. The existence of $Z$ provides that $\T(X_s,x)\neq \emptyset$. By semicontinuity we also know that $h^1(N_{X_s}(-Z))=0$.\\
Finally assume $d=(2x+4+3s)+t$, for some $t>0$. Then we obtain the required curve in $H(d,g,3)$ by induction on $t$. We start with $X_{s,0}=X_s$. Then we construct $X_{s,t+1}$ from $X_{s,t}$ by adding a line $\ell$ which meets $X_{s,t}$ at one general point, so that $X_{s,t}\cup \ell$ is nodal, and taking a smooth deformation which fixes $Z$. \\
The condition $h^1(\mathcal{O}_X(2)) =0$ follows  by applying several times Lemma \ref{eae1} applied to smooth rational curves of degree $\leq 3$. \hfil\qed

%\begin{acknowledgements}\phantom{a}\\
%Both authors are members of the Italian GNSAGA.\\
%The authors have no competing interests to declare that are relevant to the content of this article. The authors have no financial or proprietary interests in any material discussed in this article.\\
%All authors certify that they have no affiliations with or involvement in any organization or entity with any financial interest or non-financial interest in the subject matter or materials discussed in this manuscript.

%\end{acknowledgements}


\bibliographystyle{amsplain}      % mathematics and physical sciences
\bibliography{biblioLuca.bib}   % name your BibTeX data base

\end{document}